\documentclass[12pt,leqno,twoside,a4paper]{article}
\usepackage{amsfonts,amsmath,amsthm,amssymb}
\usepackage[english]{babel}
\newtheorem{theorem}{Theorem}[section]
\newtheorem{prop}[theorem]{Proposition}
\newtheorem{lemma}[theorem]{Lemma}
\newtheorem{corollary}[theorem]{Corollary}
\theoremstyle{definition}

\newtheorem{example}[theorem]{Example}
\theoremstyle{remark}

\numberwithin{equation}{section}

\newcommand{\de}{\delta}
\newcommand{\si}{\sigma}

\newcommand{\an}[2]{\mbox{\rm rann}_{#1}({#2})}

\newcommand{\ud}[1]{\mbox{\rm udim}{#1}}

\begin{document}
\title{\bf GOLDIE CONDITIONS FOR  ORE EXTENSIONS OVER SEMIPRIME RINGS}
\author{ {\bf  Andr\'{e} Leroy$^\dagger$ and Jerzy
Matczuk$^\ddagger$ }\vspace{12pt}
\\   $^\dagger$ Universit\'{e} d'Artois,  Facult\'{e} Jean Perrin\\
Rue Jean Souvraz  62307 Lens, France\\
     E-mail: leroy@euler.univ-artois.fr\\
\vspace{6pt}
\\  $^\ddagger$Institute of Mathematics, Warsaw University,\\
 Banacha 2, 02-097 Warsaw, Poland\\
  E-mail: jmatczuk@mimuw.edu.pl}
\date{ }
\maketitle\markboth{ \bf A. LEROY AND J.MATCZUK}{ \bf GOLDIE
CONDITIONS FOR  ORE EXTENSIONS}


\begin{abstract}
 Let $R$ be a ring, $\sigma$ an injective endomorphism
 of $R$ and $\de$ a $\si$-derivation of $R$.
We prove that if $R$ is semiprime left Goldie then the same holds
for the Ore extension $R[x;\si,\de]$ and both rings have the same
left uniform dimension.
\end{abstract}

\section{Introduction}

Throughout the paper $R$ will always denote an associative ring
with unity. $R[x;\si,\de]$ will stand for the Ore extension of
$R$, where $\si$ is an injective endomorphism and $\de$ a
$\si$-derivation of $R$.

 It is well-known that when $R$ is   semiprime
and $\si$ is an automorphism then the Ore extension $R[x;\si;\de]$
is a semiprime left (right) Goldie ring if and only if the ring
$R$ is such and then $R$ and $R[x;\si,\de]$ have the same uniform
dimensions. On the other hand, the easiest examples of left but
not right Ore domains can be constructed  as Ore extensions of the
form $R[x;\si]$, where $R$ is a field and $\si$ is an injective
endomorphism of $R$ which is not onto. In such a case left uniform
dimensions of $R$ and $R[x;\si]$ are equal to one but right
uniform dimensions of those rings are equal to one and infinity,
respectively.

The aim of the paper is to show that when $R$ is a semiprime left
Goldie ring and  $\si$ is injective   then the Ore extension
$R[x;\si,\de]$ is also   semiprime left Goldie  and both rings
have the same left uniform dimensions. Contrary to the
automorphism case, the fact that $R[x;\si,\de]$ is semiprime left
Goldie and $R$ is semiprime  does not imply that $R$ is left
Goldie.

 In section 2 below, with the help of certain
 classification results concerning
injective endomorphisms and skew derivations, obtained by Cauchon
and Robson (Cf. \cite{CR}), we study Jordan extensions.  In
particular, the Jordan extensions of semisimple rings are
completely described. These extensions are important tools for
proving our main results in section 3.  In this section we also
show that the ring $R$ and its Jordan extension $A$ have the same
left uniform dimension provided $R$ is semiprime left Goldie.

\section{Jordan Extensions}
Let  $R$ be ring with a fixed     injective endomorphism $\si$. We
say that an over-ring $A$ of $R$ is a Jordan extension of $R$ if
$\si$ extends to an automorphism of $A$ and $A=
\bigcup_{i=0}^\infty \sigma^{-i}(R)$. In this case we will write
$R\subseteq_\si A$.

 Jordan showed (Cf.\cite{Jo}), with the use of  left
localization of the Ore extension $R[x;\sigma]$ with respect to
the set of powers of $x$, that for any pair $(R,\si)$,  such an
extension $A$  always exists. Then he  studied the passage of
various algebraic properties from $R$ to $A$.

It is easy to observe that the Jordan extension $A$ of $R$ is an
universal object, i.e. if $R \subseteq_\si A$ and $R\subseteq_\si
A'$ are two Jordan extensions of $R$ then  the rings $A$ and $A'$
are isomorphic, by an isomorphism which is identity on $R$.

In this section we will analyse the Jordan extension and collect
some basic facts preparing the ground for  results contained in
the next section.

Throughout the paper $\mathbb{N}$ will denote the set of all
natural numbers and $\mathbb{N}_0$ will stand for the set of all
nonnegative integers.

If $\si$ is an endomorphism of the ring $R$, then $\si $ naturally
induces the  endomorphism $M_n(\si)$ of the full $n\times n$
matrix ring $M_n(R)$.

For an invertible element $u\in R$, $I_u$ will denote the inner
automorphism of $R$ adjoint to $u$, i.e. $I_u(r)=u^{-1}ru$ for all
$r\in R$.

\begin{lemma}\label{induced Jordan extensions}
Suppose that $R\subseteq _{\si}A$ is a Jordan extension and $u\in
R$ is an invertible element. Then:
\begin{enumerate}

  \item  For any $n\in \mathbb{N}$,   $M_n(R)\subseteq _{M_n(\si)}M_n(A)$ is a Jordan
extension.
\item If $\tau$ is an automorphism of $A$ such that
$\si\tau=\tau\si$, then $\tau(R)\subseteq _{\si}A$ is a Jordan
extension. If additionally $\tau(R)\subseteq R$ then  $R\subseteq
_{\si\tau}A$  is also a Jordan extension.
  \item $R\subseteq _{ I_u\si} A$ is a Jordan extension.
\end{enumerate}
\end{lemma}
\begin{proof}
 The statements (1) (2) are clear.

 (3) Obviously $I_u\si$ is an automorphism of $A$ such that $I_u\si(R)\subseteq
 R$.

 Let $a\in A$. Then there exists $n\in \mathbb{N}$ such that $\si^n(a)\in
 R$ and, consequently,  $(I_u\si)^n(a)= u^{-1}\si(u^{-1})\ldots
 \si^{n-1}(u^{-1})\si^n(a)\si^{n-1}(u)\ldots \si(u)u\in R$ follows.
\end{proof}

\begin{prop}
\label{different jordan extensions} Let $A$ be an over-ring of a
ring  $R$ and $\si$ an injective endomorphism of $R$.
 The following statements are equivalent:
\begin{enumerate}
\item For all $k\in \mathbb{N}_0$ and $n\in \mathbb{N}$,  $  \si^k(R) \subseteq_{\si^n}
A$ is a Jordan extension.
\item For all $n \in \mathbb N \, , \; R \subseteq_{\sigma^n}
A$ is a Jordan extension.
\item $R\subseteq_{\sigma} A$ is a Jordan extension.
\item There exists $n_0 \in \mathbb N$ such that $R
\subseteq_{\sigma^{n_0}} A $ is a Jordan extension.
\item There exist $k_0\in\mathbb{N}_0$ and $n_0 \in \mathbb{N}$ such that $\sigma^{k_0}(R)
\subseteq_{\sigma^{n_0}} A$ is a Jordan extension.
\end{enumerate}

\end{prop}
\begin{proof}
We only need to prove the implication $(5) \Rightarrow (1)$.
Suppose that $k_0\in\mathbb{N}_0$ and $n_0 \in \mathbb{N}$ are
such that  $\sigma^{k_0}(R) \subseteq_{\sigma^{n_0}} A$ is a
Jordan extension. In view of Lemma \ref{induced Jordan
extensions}(2), it is enough to show that $R\subseteq_{\si}A$ is a
Jordan extension. For doing so, let us extend the injective
endomorphism $\si$ of $R$ to an automorphism $\bar{\si}$ of $A$.
By assumption, we know that for any $a\in A$  there exists
$l(a)\in \mathbb N$ such that $b:=(\sigma^{n_0})^{l(a)}(a) \in
\sigma^{k_0}(R)$. Then we set $\overline{\sigma}(a):=
(\sigma^{n_0})^{-l}(\sigma(b))\in A$. This definition makes sense,
i.e. it does not depend on the choice of $l(a) $, since
$\sigma^{n_0}$ is an automorphism of $A$.

It is easy to check that $\bar{\si}$ is an endomorphism of $A$
such that for any $r\in R\;\;$ $\overline{\sigma}(r)=\sigma(r) $
and $\sigma^{-n_0}(\sigma (r))=
\overline{\sigma}(\sigma^{-n_0}(r))$. Using this, one can check
that $(\overline{\sigma})^{n_0}=(\sigma)^{n_0}$. This yields that
$\bar{\si}$ is an automorphism of $A$ such that
$\si=\bar{\si}|_R$. Then it is clear that $R\subseteq _{\si}A$ is
a Jordan extension.

\end{proof}

\begin{prop}
\label{JC extension of a product}
 Let $R$ be a
ring with an injective endomorphism $\si$.  Then:

\begin{enumerate}
\item If $R$ is a simple artinian ring and $R\subseteq_\sigma A$
is the corresponding Jordan extension, then $A$ is also simple
artinian.
\item Suppose that there is a finite bound on the cardinality of sets of orthogonal
 idempotents of $R$. Thus  $R=\prod_{i=1}^{n}R_i$ is a finite
 product of indecomposable rings. Then:
\begin{enumerate}

\item There exists a permutation $\rho\in \mathcal{S}_n$ of the index set
$ \{1,\dots,n\}$ such that $\sigma (R_i)\subseteq R_{\rho(i)}$.
\item There exists $l\in \mathbb N$ such that for any $1\leq i\leq
n$, the restriction $\si_i^l$ of $\si^l$ to $R_i$ is an injective
endomorphism of $R_i$.
\item Suppose that $l\in \mathbb N$ is as in  (b) above  and
for any $1\leq i\leq n$, $ R_i \subseteq_{\si^l_i} A_i$ is a
Jordan extension.  Then $\sigma$ can be extended to an
automorphism of  $\;\prod_{i=1}^n A_i$ and $\; \prod_{i=1}^n R_i
\subseteq_{\sigma} \prod_{i=1}^n A_i$ is a Jordan extension.
\end{enumerate}
\end{enumerate}
\end{prop}
\begin{proof}

(1) Suppose that $R$ is simple artinian and $R\subseteq_{\si}A$ is
the corresponding Jordan extension. If $R=D$ is a division ring,
then it is known and easy to show that $A$ is also a division ring
in this case.  Let  $R=M_n(D)$ where $D$ is a division ring. Then,
by Theorem 2.4 of \cite{CR},  there exists an endomorphism $\tau$
of $D$ and an invertible element $u\in R$  such that $\sigma =
I_uM_n(\tau)$. Let $K$ be the division ring such that $D
\subseteq_\tau K$ is a Jordan extension. Then, by Lemma
\ref{induced Jordan extensions}(1) and (3),  we easily obtain that
$M_n(D)\subseteq_\sigma M_n(K)$ is a Jordan extension. This means
that $A=M_n(K)$ is simple artinian.

(2) Suppose that $R=\prod_{i=1}^{n}R_i$. \\The statement (a) is
exactly Lemma 1.1 from \cite{CR}.

 (b) Let $l$ denote the order of $\rho $ in $ \mathcal S_n$.
  By (a) above $\si^l(R_i)\subseteq R_i$ for any $1\leq i\leq n$ and clearly  the
  restriction $\si^l_i$ is monic. Now, by Theorem 1.3\cite{CR}
  applied to $\si^l$, we obtain $\si^l_i(e_i)=e_i$, where $e_i$
  denotes  the unity of $R_i$.

(c) Let $1\leq i\leq n$ and $l$ be as in (b). In view of (b), we
can consider the Jordan extensions $R_i \subseteq_{\sigma^l_i} A_i
$ for $1\leq i\leq n$. Then $\prod^n_{i=1} R_i \subseteq_{\si^l}
\prod^n_{i=1} A_i$ is also a Jordan extension. Now, Proposition
\ref{different jordan extensions} completes the proof.
\end{proof}

The above Proposition gives us immediately the following:

\begin{corollary}\label{Cor. semisimple J ext}
Let a Jordan extension $R\subseteq_{\si} A$ of a semisimple ring
$R$ be given. Then
  $A$ is   semisimple. In fact, if
 $R=\prod_{i=1}^{n}M_{n_i}(D_i)$ for  some division rings $D_i$, then
 $A=\prod_{i=1}^{n}M_{n_i}(K_i)$ for suitable division rings $K_i$.
 \end{corollary}

For a  semiprime left Goldie ring $R$, $Q(R)$ will denote the
classical left quotient ring of $R$. Recall that, by Goldie's
Theorem, $Q(R)$ is a semisimple ring.

The following results will be useful for our purposes.

\begin{prop}\label{from Jordan}
Let $R$ be a semiprime left Goldie ring with an injective endomorphism
$\si$ and $R \subseteq_{\si} A$ be the corresponding Jordan
extension. Then:
\begin{enumerate}
\item $\si(\mathcal{C} )\subseteq \mathcal{C}$, where
$\mathcal{C}$ denotes the set of all regular elements of $R$.
\item $\si$ can be uniquely extended to an injective endomorphism of
$Q(R)$.
\item $A$ is a semiprime left Goldie ring and $Q(R)\subseteq_\sigma Q(A)$
is a Jordan extension.
\item Every $\si$-derivation $\de$   of $R$ has a
unique extension to a $\si$-derivation of $Q(R)$.

\end{enumerate}
\end{prop}
\begin{proof}
    The first statement is a special case of a result of
Jategaonkar (Cf. Proposition 2.4 \cite{Ja}) which states that
$\si( \mathcal{C})\subseteq\mathcal{C} $ when $R$ has left
artinian quotient ring.
  (2), (4). The fact that $\si$ can be extended to $Q(R)$
is a well-known consequence of (1). Then, it is also known that
every $\si$-derivation $\de$ of $R$ extends uniquely to a
$\si$-derivation of $Q(R)$ and
$\de(c^{-1})=-\si(c^{-1})\de(c)c^{-1}$ for all $c\in\mathcal{C}$.

 (3).  We know that $Q(R)$ is semisimple.  Let
$Q(R) \subseteq_\si A(Q(R))$ be a Jordan extension for $\si$
extended to $Q(R)$. Then, by Corollary \ref{Cor. semisimple J
ext},  $A(Q(R))$ is also semisimple.

  For any $x\in A(Q(R))$
there exist $n \in \mathbb N, \; c\in \mathcal C$ and $r\in R$
such that $\si^n(x)= c^{-1}r$, i.e. $x =
\si^{-n}(c^{-1})\si^{-n}(r)$. This shows that
$A:=\bigcup_{i=0}^\infty \si^{-i}(R) \subseteq A(Q(R))$ is a left
order in the semisimple ring $A(Q(R))$.  Now,  Theorem 3.1.7 in
\cite{MR} yields  that $A$ is semiprime left Goldie. Then
$A(Q(R))=Q(A)$ easily follows.
\end{proof}

The statement (3) from  the above proposition was also obtained,
using other methods, by Jordan in \cite{Jo} (see \cite{Jo}
Corollary 7.5, Proposition 7.1 and Theorem 7.2).  The above proof
was given both for completeness of the presentation and as an
application of Corollary \ref{Cor. semisimple J ext}.

The following proposition shows that in case the $\si$-derivation
$\de$ is $q$-quant\-ized, i.e. $\de\si=q\si\de$ for some central,
$\si$ and $\de$ invariant element $q\in R$, a Jordan extension $R
\subseteq_\si A$ leads to a Jordan extension $R[x;\si,\de]
\subseteq_\si A[x;\si,\de]$.

\begin{prop}
\label{R and A Q(R) and Q(A)} Let $R \subseteq_{\sigma} A$ be a
Jordan extension and $\delta$ be a $q$-quantized
$\sigma$-derivation of $R$. Then:
\begin{enumerate}
\item $\delta$ can be uniquely extended to a $q$-quantized $\sigma$ derivation of $A$.
 \item $\sigma$ can be extended to an injective endomorphism of $R[x;\si,\de] $.
 Moreover,
$R[x;\sigma,\delta] \subseteq_\sigma A[x;\sigma,\delta]$ is a
Jordan extension.
\end{enumerate}
\end{prop}
\begin{proof}
(1).  Suppose that $\bar{\de}$ is an extension of $\de$ to a
$q$-quantized $\si$-derivation of $A$. Let $a \in A$ and $n\in
\mathbb{N}_0$ be such  such that $\sigma^n(a)\in R$. Then
$\bar{\de} (a)=
q^{-n}\sigma^{-n}(\bar{\de}(\sigma^n(a)))=q^{-n}\sigma^{-n}(\de(\sigma^n(a)))$.
This shows that $\bar{\de}$ is uniquely determined by $\de$  and $
\si$.

Notice also that if $a\in A$ and $n,m\in\mathbb{N}_0$ are such
that $\si^n(a),\si^m(a)\in R$, then
$q^{-n}\sigma^{-n}(\de(\sigma^n(a)))=q^{-m}\sigma^{-m}(\de(\sigma^m(a)))$.
Now, it is standard to check that  $\de \colon A\rightarrow A$
given by $ \de(a) = q^{-n}\sigma^{-n}(\delta(\sigma^n(a)))$, where
$n\in \mathbb{N}_0$ is such that $\si^n(a)\in R$ is a well defined
$q$-quantized $\si$-derivation of $A$.

(2). The fact that $\sigma$ can be extended to an injective
endomorphism of $R[x;\si, \de]$ is part of folklore: just define
$\sigma(x)=q^{-1}x$.

Let $p=\sum_{i=0}^la_ix^i \in A[x;\sigma,\delta]$. By assumption,
 there are   $n_i\in \mathbb{N}_0$, with $0\leq i\leq l$, such that
$\sigma^{n_i}(a_i)\in R$. Then $\sigma^n(p)\in R[x;\sigma,\delta]$
where $n=\max \{n_i \mid 0\leq i\leq l\}$.   This easily yields
 that $R[x;\sigma,\delta] \subseteq_\sigma
A[x;\sigma,\delta]$ is a Jordan extension.
\end{proof}

\section{Main Results}
We begin this section with a description of skew polynomial rings
over semisimple rings.
\begin{lemma}\label{decomposition of skew pol ring}
 Let $R$ be a semisimple ring, $\si$ and $\de$ an injective
 endomorphism and a $\si$-derivation of $R$, respectively. Then
 either
\begin{enumerate}
  \item   $R$ is simple artinian and   there exists a division ring $D$ with
  an  endomorphism $\si'$ and $\si$-derivation $\de'$ such that
   $R[x,\si,\de]\simeq M_m(D[y,\si',\de'])$ for some $m\in
   \mathbb{N}$.
   \end{enumerate}
   or
\begin{enumerate}
  \item[2]    There exists a  ring decomposition
 $R[x;\si,\de]=\prod^k_{j=0}B_j[x_j,\si_j,\de_j]$ such that:
\begin{enumerate}
  \item  $\si_j$ is an injective endomorphism of $B_j$ for any $1\leq j\leq k$.
  \item If for some $1\leq j\leq k $ $\;B_j$ is not a simple ring,
  then $\de_j=0$.
 \end{enumerate}
\end{enumerate}
\end{lemma}
\begin{proof} (1) Suppose $R$ is simple artinian. Then $R=M_m(D)$ for
some division ring $D$ and $m\in\mathbb{N}$ and the statement (1)
is a particular case of Theorem 3.2  from \cite{CR}.

(2) Let $1=e_1+ \ldots +e_n$ be the decomposition of $1$ into the
sum of    central primitive orthogonal idempotents of $R$. Then
$R=\prod_{i=1}^nRe_i$ and each $Re_i$ is simple artinian.  By
Proposition \ref{JC extension of a product}(2)(a), $\si$ induces a
permutation $\rho$ of the index set
 $\{1,\ldots ,n\}$. Let $\mathcal{O}_1,\ldots
,\mathcal{O}_k$ denote the orbits of this action and set
$B_j=\prod_{i\in \mathcal{O}_j}Re_i$ for $1\leq j\leq k$. Let
$\si_j$, $\de_j$ be the restriction of $\si$ and $\de$ to $B_j$.
Then, by Theorem 1.3 from \cite{CR}, $\si_j$ is an injective
endomorphism of $B_j$ and $\de_j$ is a $\si_j$-derivation of
$B_j$. Therefore, we can decompose the ring $R[x;\si,\de]$ in the
following way: $$R[x;\si ,\de]=\prod_{j=1}^k(B_j[y_j;\si_j
,\de_j]).$$ When the cardinality $\#\mathcal{O}_j$ of the orbit
$\mathcal{O}_j$ is equal to 1, then $B_j$ is simple and we can set
$x_j=y_j$ in this case. When $\#\mathcal{O}_j>1$  then, by Lemma
1.4 \cite{CR}, $\de_j$ is an inner $\si_j$-derivation of $B_j$,
i.e. there is $b\in B_j$ such that $\de(r)=br-\si(r)b$ for any
$r\in B_j$. Then   $B_j[y_j;\si_j,\de_j]$ is isomorphic to the
ring $B_j[x_j;\si_j]$, where $x_j=y_j-b$. This completes the proof
of the lemma.
\end{proof}

In the sequel, the left uniform dimension of a ring $R$ is denoted
by $\ud{R}$. The following result was obtained by Mushrub in
\cite{Mu}. The original argument was lengthy, thus we present a
new very short  proof.

\begin{lemma}
\label{udim when delta = 0} Let $R$ be a ring and $R\subseteq_\si
A$ a Jordan extension associated to an injective endomorphism of
$R$.  Then $\ud{R[x;\si]}=\ud{A[x;\si]}=\ud{A}$.
\end{lemma}

\begin{proof} Let $S=\{x^n\mid n\geq 0\}$. It is known and easy
to check (Cf. \cite{Jo}) that $S$ is a left Ore set in $R[x;\si]$
and $S^{-1}R[x;\si]$ is isomorphic to
$A[x,x^{-1};\si]=S^{-1}A[x;\si]$. The left uniform dimension is
preserved under left localizations with respect to Ore sets of
regular elements (Cf. Lemma 2.2.12 \cite{MR}). Hence
$\ud{R[x;\si]}=\ud{A[x;\si]}$ follows.

Now,  as $\si$ is an automorphism of $A$,  a classical result of
Shock (Cf. \cite{Sh}, \cite{Ma}) says that $\ud{A}=\ud{A[x;\si]}$.
\end{proof}

\begin{lemma}
\label{udim of semi simple}
 Let $R$ be a semisimple  ring, $\si$ an injective endomorphism  and
 $\delta$ a $\sigma$-derivation of $R$.
Then $\ud{R[x;\si,\de]}=\ud{R}$.

\end{lemma}
\begin{proof}
Recall that if a ring $B$ is isomorphic to $\prod_{j=1}^kB_j$,
then $\;\ud{B}=\sum_{j=1}^k\ud{B_j}$. Hence, in virtue of Lemma
\ref{decomposition of skew pol ring}, it is enough to prove the
lemma in two special cases:  when $\de=0$ and  when $R$ is simple
artinian.

{\em Case 1}. Suppose that $\de=0$, i.e. $R[x;\si,\de]=R[x;\si]$.
Let $R\subseteq_\sigma A$ be the Jordan extension.  Since $R$ is
semisimple,  Corollary \ref{Cor. semisimple J ext} implies that
$\ud{A}=\ud{R}$.  The thesis is now clear, thanks to Lemma
\ref{udim when delta = 0}.

{\em Case 2}. Suppose that $R$ is simple artinian, i.e $R\simeq
M_m(D)$ for some division ring $D$. Thus, by Lemma
\ref{decomposition of skew pol ring} the ring $R[x;\si,\de]$ can
be presented in the form $ M_m(D[y;\si',\de'])$ . Since $D$ is a
division ring, $D[y;\si',\de']$ is a principal left ideal domain,
so $\ud{D[y;\si',\de']}=1$ and $\ud{R[x;\si,\de]}=m=\ud{R}$
follows. This completes the proof.
\end{proof}

\begin{theorem}\label{udim}
Let $R$ be a semiprime left Goldie ring with   an injective
endomorphism $\si$ and $R \subseteq_{\sigma} A$ be the
corresponding Jordan extension.  Then for any $\si$-derivation
$\de$  of $R$  $\;\ud{R[x;\si,\de]}=\ud{R}=\ud{A}$.
\end{theorem}
\begin{proof} By Proposition \ref{from Jordan}, $\si$ extends to
an injective endomorphism of $Q(R)$ and $\de$ extends to a
$\si$-derivation of $Q(R)$. Thus, we can consider the ring
extension $R[x;\si ,\de]\subseteq Q(R)[x;\si,\de]$.

Clearly all elements from the set  $\mathcal{C}$  of all regular
elements of $R$ are invertible in $Q(R)[x;\si,\de]$ and every
element from $Q(R)[x;\si ,\de]$ can be presented in the form
$c^{-1}p$ for some $c\in \mathcal{C}$ and $p\in R[x;\si ,\de]$.
This means that $\mathcal{C}$ is a left Ore set in $R[x;\si ,\de]$
and $\mathcal{C}^{-1}(R[x;\si ,\de])=Q(R)[x;\si,\de]$. Now, by
\cite{MR} Lemma 2.2.12, we obtain:
$$\ud{R[x;\si,\de]}=\ud{Q(R)[x;\si,\de]}\;\; \mbox{and}   \;\;
\ud{R}=\ud{Q(R)}.$$ By Proposition \ref{from Jordan} we know that
$A$ is semiprime left Goldie, so we also have $\ud{A}=\ud{Q(A)}$.

Notice  that the same proposition yields that
$Q(R)\subseteq_{\si}Q(A)$ is a Jordan extension. Since both $Q(R)$
and $Q(A)$ are semisimple rings, Corollary \ref{Cor. semisimple J
ext} and Lemma \ref{udim of semi simple} give us
$$\ud{Q(R)}=\ud{Q(A)} \;\; \mbox{and} \;\;\ud{Q(R)[x;\si,\de]}=\ud{Q(R)},$$ respectively.
This implies the thesis.
\end{proof}

Mushrub in \cite{Mu} investigated the left uniform dimension of
skew polynomial rings $R[x;\si]$, where $\si$ is an injective
endomorphism of $R$. He proved, in particular, that when $R$ is a
left Ore domain (i.e. a domain with $\ud{R}=1$), then $R[x;\si]$
is   a left Ore domain. He also constructed a series of examples
showing that:

\begin{enumerate}
  \item   For any $n\in\mathbb{N}$ there is a
commutative ring $R$ (not semiprime) with an injective
endomorphism $\si$, such that $\ud{R}=n$ and $\ud{R[x;\si]}=1$.
  \item There exists a domain $R$ with an injective
endomorphism $\si$ such that $\ud{R[x;\si]=1}$ but $R$ has
infinite both left and right uniform dimensions.
 \end{enumerate}

 He posed a
question  whether $\ud{R}=\ud{R[x;\si]}$ provided $R$ is a
semi\-prime ring of finite left Goldie dimension. The above
theorem gives a positive answer to this question for $R$
satisfying the ACC on left annihilators. We also have the
following:
\begin{corollary}
 Let $R$ be a left Ore domain with an injective endomorphism
 $\si$. Then $R[x;\si,\de]$ is a left Ore domain, for any $\si$-derivation $\de$ of
 $R$.
\end{corollary}
\begin{proof}
 Since $\si$ is injective, $R[x;\si,\de]$ is a domain. Now the
 thesis is a direct consequence of Theorem \ref{udim}
\end{proof}

Using Lemma \ref{decomposition of skew pol ring} one can easily
show that when $R$ is a semiprime left Goldie ring  then the skew
polynomial ring $R[x;\si,\de]$ is semiprime provided $\si$ in
injective. The following lemma is slightly more general.
\begin{prop}
\label{semiprimeness} Suppose that the  ring $R$   satisfies the
{\rm ACC} on left annihilators.  Let $\si$ and $ \de$ stand for an
injective endomorphism and a $\si$-derivation of $R$,
respectively.  Then the Ore extension $R[x;\si,\de]$ is prime
(semiprime), provided $R$ is prime (semiprime).
\end{prop}
\begin{proof}
Let $I$ be an ideal of $R[x;\si,\de]$.  For any $n\geq 0 $ define
 $ I_n=\{a\in R\mid$ either
$a=0$ or $a$ is the leading coefficient of some polynomial from
$I$ of degree $n\}$.  Clearly $\{I_n\}_{n\geq 0}$ is an ascending
sequence of left ideals of $R$ such that $\sigma^l(I_n) \subseteq
I_{n+l}$ for any $n,l\geq 0$. Since $R$ satisfies the ACC on left
annihilators, it satisfies the DCC on right annihilators.
Therefore, there exists $n_0\ge 0$ such that $
\an{R}{I_n}=\an{R}{I_{2n}}$ for any $n\geq n_0$.

 Suppose that $R$ is prime and let $J$ be a
nonzero ideal of $R[x;\si,\de]$. By the considerations above,
there is $m\geq 0$ such that for any $n\geq m$  $
\an{R}{I_n}=\an{R}{I_{2n}}$ and $J_n\ne 0$. Assume $IJ=0$. Then
$I_n\si^n(J_n)=0$. This means that $\si^n(J_n) \subseteq
\an{R}{I_n}=\an{R}{I_{2n}}$ and $I_{2n}\si^n(J_n)=0$ follows.
Since $\si^n(I_n) \subseteq I_{2n}$, we get
$\si^n(I_nJ_n)=\si^n(I_n)\si^n(J_n)=0$. This leads to $I_nJ_n = 0$
for any $n\geq m$, as $\si$ is injective. Since $J_n\ne 0$,
primeness of $R$ yields $I_n$=0 for all $n\geq m$ and $I=0$
follows. This shows that $R[x;\si,\de]$ is prime.

The same argument applied to $J=I$ shows that $R[x;\si,\de]$ is
semiprime provided $R$ is semiprime.
\end{proof}
In the case $\si$ is an automorphism then the Ore extension
$R[x;\si,\de]$ is always prime when $R$ is prime. However
semiprimeness of $R$ does not imply semiprimeness of
$R[x;\si,\de]$. In the case $\si$ is an injective endomorphism,
the situation is more complex, as the following example shows.
Namely, there exists a prime ring $R$ such that $R[x;\si,\de]$ is
not semiprime.

\begin{example} Let $R$ be a subset of $\mathbb{N}\times
\mathbb{N}$ matrices over a field $K$ defined as follows
 $R=\{M\mid M=\sum_{i,j=1}^na_{ij}e_{ij}+
a\sum_{i=n+1}^{\infty}e_{ii}$ for some $n\in\mathbb{N}$ and
$a_{ij}, a\in K\}$, where $\{e_{ij}\}_{i,j\in\mathbb{N}}$ denotes
the set of matrix units. Then $R$ is a prime unital ring and it is
easy to check that the map $\si\colon R\rightarrow R$ given by
$$\si(\sum_{i,j=1}^na_{ij}e_{ij}+
a\sum_{i=n+1}^{\infty}e_{ii})=ae_{11}+
\sum_{i,j=1}^na_{ij}e_{i+1,j+1}+a\sum_{i=n+2}^{\infty}e_{ii}$$ is
an injective endomorphism of $R$.

Notice that $e_{11}\si (R)=Ke_{11}$.  Therefore, for any $k\geq
0$, we have in $R[x;\si]$ :
$$e_{11}xRx^ke_{11}=Ke_{11}e_{2+k,2+k}x^{k+1}=0.$$
This shows that $e_{11}xR[x;\si]$ is a nilpotent left ideal of
$R[x;\si]$, i.e.  $R[x;\si]$ is not semiprime.
\end{example}

Let $R\subseteq_{\si}A$  be the Jordan extension where $R$ and
$\si$ are as in the above example. Then, by Proposition \ref{R and
A Q(R) and Q(A)}, we can consider the Jordan extension
$R[x;\si]\subseteq_{\si} A[x,\si]$. Then, as we have seen,
$R[x;\si]$ is not semiprime but $A[x,\si]$ is a prime ring since
$A$ is prime and $\si$ is an automorphism of $A$.

 We close the paper with the following:
\begin{theorem}\label{goldie}
 Let $R$ be a semiprime left Goldie ring, $\si$, $\de$ an
injective endomorphism and a $\si$-derivation of $R$,
respectively.  Then $R[x;\si,\de]$ is also a semiprime left Goldie
ring.
\end{theorem}
\begin{proof}
 By Proposition \ref{semiprimeness} and Theorem \ref{udim} we know
 that $R[x;\si,\de]$ is a semi\-prime ring of finite left Goldie
 dimension. Thus, in order to complete the proof, it is enough to
 show that $R[x;\si,\de]$ satisfies the ACC on left annihilators.
 As we have seen in Proposition \ref{from Jordan}, $\si$ and $\de$
 can be extended to an injective endomorphism
 and a $\si$-derivation of $Q(R)$. Since the ACC on left annihilators is inherited
 on subrings, it is enough to prove that $Q(R)[x;\si,\de]$ satisfies
 the ACC on left annihilators. This means that without loss of generality
  we may assume that $R$ is semisimple. Then, by making use of
  Lemma \ref{decomposition of skew pol ring}, it is enough to
  consider only two cases: when $R$ is simple
artinian and when $\de=0$.

{\em Case 1}. Suppose $R$ is simple artinian. Then, by Lemma
\ref{decomposition of skew pol ring}, $R[x;\si,\de]$ is isomorphic
to $M_m(D[y,\si',\de'])$ for some $m$, where $D$ is a division
ring. $D[y,\si',\de']$ is a left principal domain  so, in
particular, it is a prime left Goldie ring. Therefore, by
Corollary 3.1.5 from \cite{MR}, $R[x;\si,\de]$ is a prime left
Goldie ring, thus satisfies the ACC on left annihilators.

{\em Case 2}. Suppose that $\de=0$, i.e. $R[x;\si,\de]=R[x;\si]$.
Let $R\subseteq_\sigma A$ be the Jordan extension. By Corollary
\ref{Cor. semisimple J ext}, $A$ is also semisimple. Now $\si$ is
an automorphism of $A$  so, by Theorem 2.6 of \cite{Ma},
$A[x;\si]$ is a semiprime  Goldie ring and $R[x;\si]$ satisfies
the ACC on left annihilators as a subring of $A[x;\si]$.
\end{proof}

If $\si$ is an automorphism of a semiprime ring $R$, then it is
known (Cf. Theorem 2.6 \cite{Ma}) that $R$ is semiprime left
Goldie if and only if the Ore extension $R[x;\si,\de]$ is
semiprime left Goldie. The quoted earlier example of Mushrub shows
that when $\si $ is just an injective endomorphism, then the above
equivalence does not hold. However, using Theorem \ref{goldie},
one can easily show that $R[x;\si, \de]$ is semiprime left Goldie
if and only if $R$ is semiprime left Goldie provided $R$ is
semiprime and $\ud{R}$ is finite.

\end{document}